\documentclass[11pt, fceqn]{article}
\usepackage[utf8]{inputenc}
\usepackage{graphicx}
\usepackage{cite}
\usepackage{amsfonts}
\usepackage{subcaption}
\usepackage{multirow}
\usepackage{tabu} 
\usepackage{makeidx} 
\usepackage{ifpdf}
\usepackage{url}

\usepackage{amsmath}
\usepackage{amssymb}
\usepackage[sort]{natbib}
\usepackage{graphicx}
\usepackage{float}
\usepackage[left = 1.00in, right = 1.00in, top = 1.00in, bottom = 1.00in]{geometry}
\usepackage{amsthm}

\usepackage{rotating}
\usepackage{algorithm}
\usepackage{algpseudocode}
\usepackage{setspace}
\usepackage{comment}
\title{Managing Access to Care at Primary Care Clinics}

\begin{document}
	
\begin{center}
	{\Large Note on the article: Optimal choice for appointment scheduling window under patient no-show behavior}\\[12pt]
		
	\footnotesize
	\mbox{\large Sina Faridimehr}\\
	Turner Broadcasting System Inc., Atlanta, GA 30318, USA. \\
	\normalsize
\end{center}

\begin{comment}
	\author[label1]{Sina Faridimehr\corref{cor1}}
	\author[label1]{Saravanan Venkatachalam}
	\author[label1]{Ratna Babu Chinnam}
	%\author{Daniel Jornada\fnref{label1}, V.Jorge Leon\corref{cor1}\fnref{label1}, Saravanan Venkatachalam \fnref{label2}}
	\address[label1]{Department of Industrial \& Systems Engineering, Wayne State University}
	%\fntext[label1]{}
	%\fntext[label2]{}
	\cortext[cor1]{corresponding author: sina.faridimehr@wayne.edu}
\end{comment}

\begin{abstract}
	In a recently published article, Liu (2016) used M/M/1/K and M/D/1/K queuing models to determine the optimal window for patient appointment scheduling considering clinic environmental factors such as panel size, provider service rate, and delay-dependent no-show probability. In this note, we show that the application of one of the no-show functions was incorrect. The required corrections are proposed.
\end{abstract}

\section{Introduction}

\begin{figure*}[!htbp]
	\centering
	\includegraphics[scale = 0.4]{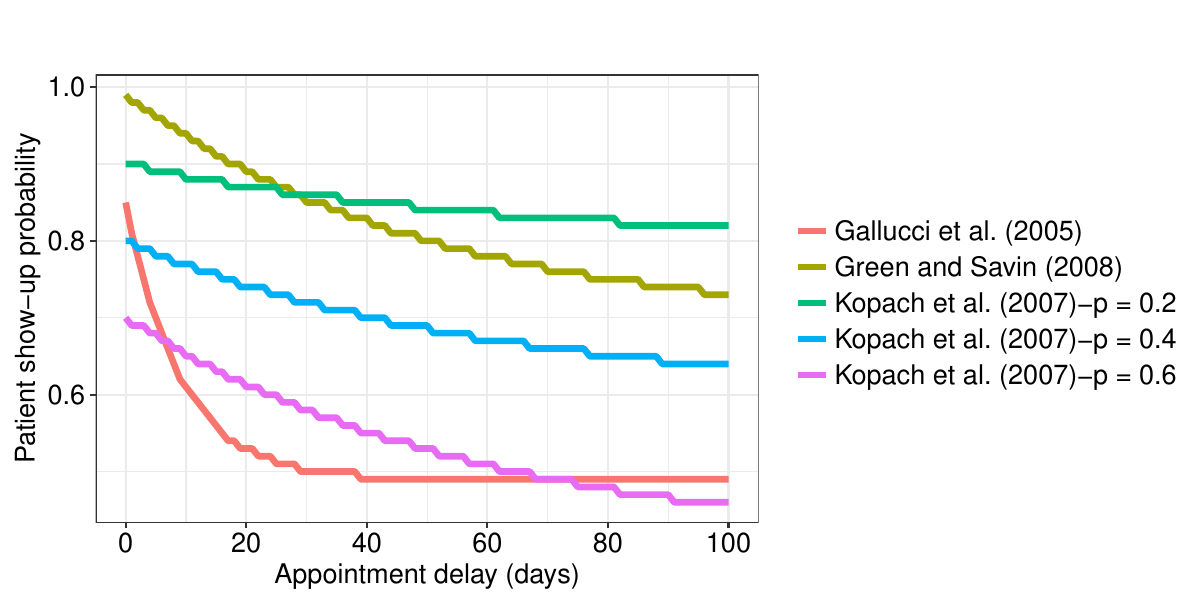}
	\caption{Patient show up probabilities}
	\label{Show}
\end{figure*}

Kopach et al. (2007) proposed the following function for the relationship between show-up rate and appointment delay:
\begin{equation}
	p_{j} = 1-p(1-0.5*e^{-0.017j})
\end{equation}
where $p$ is the patient estimated no-show rate and j is appointment lead time. However, the following function is used by Liu (2016) as show-up rate which artificially shows higher no-show rate:
\begin{equation}
p_{j} = 0.5*e^{-0.017j}
\end{equation}
We use the correct function proposed by Kopach et al. (2007) to show the effect of this no-show behavior on optimal appointment scheduling window in M/M/1/K and M/D/1/K queuing systems. In the first case, provider has no operational lever to control in order to deal with patient no-show problem whereas in the second case provider has the ability to adjust panel size and overbooking level. As noted in Liu (2016), patients in different clinical settings have different no-show probability. Hence, we use three different levels of 0.2, 0.4 and 0.6 as estimated no-show probability in no-show function proposed by Kopach et al. (2007). Figure \ref{Show} shows the sensitivity of patients to delay under three different proposed no-show function. For the sake of completeness, we are showing Liu (2016) results for the other two no-show functions.

\section{Comparison of M/M/1/K and M/D/1/K Models}

Liu (2016) proposed the following objective function in order to find the optimal value of appointment scheduling window:
\begin{equation}
T(K) = \lambda \sum_{j=0}^{K-1}\Pi_{j}(K)q{j} + \mu\xi\Pi_{0}(K) - \lambda\theta\Pi_{K}(K)
\end{equation}
where the first term is the revenue obtained by visiting patients that showed up for their appointment, second term is the revenue obtained by doing ancillary tasks when there is no patient in clinic and the third term is the cost of rejecting patient requests due to limited capacity in queue.

Table \ref{OptimalK1} and Table \ref{OptimalK2} compare the optimal value for appointment scheduling window M/M/1/K and M/D/1/K queuing systems, respectively. As in Liu (2016), we consider service rate $\mu = 20$ and compare the optimal value for scheduling window for different values of demand rate $\lambda$, ancillary task revenue rate $\xi$ and patient request rejection penalty $\theta$. We use the same notations as Liu (2016).

\begin{table*}[!htbp]
	\caption{Optimal Scheduling Window in M/M/1/K Model}
	\centering
	\begin{tabular}{|c c c c c c c|} 
		\hline
		$(\theta,\xi)$ & $\lambda$ & $K_{M,K_{0.2}}$ & $K_{M,K_{0.4}}$ & $K_{M,K_{0.6}}$ & $K_{M,G}$ & $K_{M,GS}$ \\
		\hline
		\multirow{4}{1.5em}{(0,0)} & 18 & $\infty$ & $\infty$ & $\infty$ & 60 & $\infty$ \\
		& 19 & $\infty$ & 280 & 160 & 40 & 200 \\
		& 19.9 & 180 & 100 & 80 & 40 & 80 \\
		& 19.99 & 160 & 100 & 80 & 40 & 80 \\
		\hline
		\multirow{4}{2em}{(0,0.5)} & 18 & $\infty$ & $\infty$ & $\infty$ & 60 & $\infty$ \\
		& 19 & $\infty$ & 280 & 160 & 40 & 200 \\
		& 19.9 & 180 & 100 & 80 & 40 & 80 \\
		& 19.99 & 160 & 100 & 80 & 40 & 80 \\
		\hline
		\multirow{4}{2em}{(1.5,0)} & 18 & $\infty$ & $\infty$ & $\infty$ & 200 & $\infty$ \\
		& 19 & $\infty$ & $\infty$ & $\infty$ & 100 & $\infty$ \\
		& 19.9 & 320 & 200 & 160 & 60 & 160 \\
		& 19.99 & 260 & 180 & 140 & 60 & 140 \\
		\hline
		\multirow{4}{3em}{(1.5,0.5)} & 18 & $\infty$ & $\infty$ & $\infty$ & $\infty$ & $\infty$ \\
		& 19 & $\infty$ & $\infty$ & $\infty$ & 160 & $\infty$ \\
		& 19.9 & 460 & 280 & 200 & 80 & 200 \\
		& 19.99 & 340 & 220 & 180 & 60 & 180 \\
		\hline
	\end{tabular}
	\label{OptimalK1}
\end{table*}

\begin{table*}[!htbp]
	\caption{Optimal Scheduling Window in M/D/1/K Model}
	\centering
	\begin{tabular}{|c c c c c c c|} 
		\hline
		$(\theta,\xi)$ & $\lambda$ & $K_{D,K_{0.2}}$ & $K_{D,K_{0.4}}$ & $K_{D,K_{0.6}}$ & $K_{D,G}$ & $K_{D,GS}$ \\
		\hline
		\multirow{4}{1.5em}{(0,0)} & 18 & $\infty$ & $\infty$ & $\infty$ & 60 & $\infty$ \\
		& 19 & $\infty$ & $\infty$ & 160 & 40 & 200 \\
		& 19.9 & 120 & 80 & 60 & 20 & 60 \\
		& 19.99 & 100 & 80 & 60 & 20 & 60 \\
		\hline
		\multirow{4}{2em}{(0,0.5)} & 18 & $\infty$ & $\infty$ & $\infty$ & 60 & $\infty$ \\
		& 19 & $\infty$ & $\infty$ & 160 & 40 & 200 \\
		& 19.9 & 120 & 80 & 60 & 20 & 60 \\
		& 19.99 & 100 & 80 & 60 & 20 & 60 \\
		\hline
		\multirow{4}{2em}{(1.5,0)} & 18 & $\infty$ & $\infty$ & $\infty$ & 160 & $\infty$ \\
		& 19 & $\infty$ & $\infty$ & $\infty$ & 80 & $\infty$ \\
		& 19.9 & 260 & 160 & 120 & 40 & 120 \\
		& 19.99 & 180 & 120 & 100 & 40 & 100 \\
		\hline
		\multirow{4}{3em}{(1.5,0.5)} & 18 & $\infty$ & $\infty$ & $\infty$ & $\infty$ & $\infty$ \\
		& 19 & $\infty$ & $\infty$ & $\infty$ & 160 & $\infty$ \\
		& 19.9 & 380 & 220 & 160 & 60 & 160 \\
		& 19.99 & 240 & 160 & 120 & 40 & 120 \\
		\hline
	\end{tabular}
	\label{OptimalK2}
\end{table*}

As Tables \ref{OptimalK1} and \ref{OptimalK2} show, since patients are less sensitive to delay when the estimated no-show rate is 0.2 or 0.4, optimal window under the no-show function proposed by Kopach et al. (2007) is larger than or equal to optimal windows resulted from other no-show functions, $K^{*}_{M,G} \leq K^{*}_{M,GS} \leq K^{*}_{M,K_{0.4}} \leq K^{*}_{M,K_{0.2}}$, but if the estimated no-show rate increases to 0.6, patients become more sensitive to delay under the no-show function of Kopach et al. (2007) and the order of optimal scheduling windows will become $K^{*}_{M,G} \leq K^{*}_{M,K_{0.6}} \leq K^{*}_{M,GS}$. We also see this pattern in M/D/1/K queuing model without any exception.

As suggested by Liu (2016), we use efficiency gain to evaluate the closeness of M/M/1/K approximations to M/D/1/K queuing system. Efficiency gain measures the loss in clinic reward by choosing optimal scheduling window resulted from M/M/1/K system rather the one from M/D/1/K system:
\begin{equation}
\Delta E = 100\% \ast \frac{T(K^*_{D})-T(K^*_{M})}{K^*_{D}}
\end{equation}
According to Table \ref{EfficiencyOptimal1} and Table \ref{EfficiencyOptimal2}, 44\% of optimal windows matches between M/M/1/K and M/D/1/K queuing systems and the average amount of efficiency loss in the rest of cases is 0.84\% which confirms that analytical results obtained from M/M/1/K model is a close approximation of M/D/1/K model.

\section{Efficiency Gains Resulted from Adopting $K^*$ without Other Operational Levers}

In this case, clinic manager has no operational lever other than appointment scheduling window to maximize long-run average net reward. We calculate $T(K)$ in the case that patients can be scheduled any time in the future and compare it to the case that the practice adopts optimal window. The efficiency gain can be defined as the percentage improvement in long-run average net reward when we use optimal scheduling window:
\begin{equation}
\Delta E = 100\% \ast \frac{T(K^*)-T(\infty)}{T(\infty)}
\end{equation} 
Tables \ref{EfficiencyOptimal1} and \ref{EfficiencyOptimal2} represent the efficiency gain by adopting the optimal scheduling window when other operational levers are not in control of clinic manager.

\begin{table*}[!htbp]
	\caption{Efficiency Gains (\%) Resulted from Adopting $K^*$ without Other Operational Levers in M/M/1/K Model}
	\centering
	\begin{tabular}{|c c c c c c c|} 
		\hline
		$(\theta,\xi)$ & $\lambda$ & $\Delta E_{M,K_{0.2}}$ & $\Delta E_{M,K_{0.4}}$ & $\Delta E_{M,K_{0.6}}$ & $\Delta E_{M,G}$ & $\Delta E_{M,GS}$ \\
		\hline
		\multirow{4}{1.5em}{(0,0)} & 18 & 0.00 & 0.00 & 0.00 & 0.00 & 0.00 \\
		& 19 & 0.00 & 0.00 & 0.00 & 0.46 & 0.00 \\
		& 19.9 & 0.56 & 1.93 & 4.02 & 21.19 & 3.02 \\
		& 19.99 & 2.23 & 6.18 &	12.07 &	42.50 & 9.08 \\
		\hline
		\multirow{4}{2em}{(0,0.5)} & 18 & 0.00 & 0.00 & 0.00 & 0.00	& 0.00 \\
		& 19 & 0.00 & 0.00 & 0.00 & 0.20 & 0.00 \\
		& 19.9 & 0.26 & 0.84 & 1.59 & 8.49 & 1.46 \\
		& 19.99 & 1.04 & 2.63 &	4.57 & 15.41 & 4.27 \\
		\hline
		\multirow{4}{2em}{(1.5,0)} & 18 & 0.00 & 0.00 & 0.00 & 0.00 & 0.00 \\
		& 19 & 0.00 & 0.00 & 0.00 & 0.02 & 0.00 \\
		& 19.9 & 0.22 & 1.09 & 2.57 & 16.67 & 2.05 \\
		& 19.99 & 1.57 & 4.91 &	10.03 & 36.67 &	7.71 \\
		\hline
		\multirow{4}{3em}{(1.5,0.5)} & 18 & 0.00 & 0.00 & 0.00 & 0.00 & 0.00 \\
		& 19 & 0.00 & 0.00 & 0.00 & 0.00 & 0.00 \\
		& 19.9 & 0.04 & 0.29 & 0.72 & 5.48 & 0.73 \\
		& 19.99 & 0.54 & 1.74 &	3.31 & 11.82 & 3.20 \\
		\hline
	\end{tabular}
	\label{EfficiencyOptimal1}
\end{table*}

\begin{table*}[!htbp]
	\caption{Efficiency Gains (\%) Resulted from Adopting $K^*$ without Other Operational Levers in M/D/1/K Model}
	\centering
	\begin{tabular}{|c c c c c c c|} 
		\hline
		$(\theta,\xi)$ & $\lambda$ & $\Delta E_{D,K_{0.2}}$ & $\Delta E_{D,K_{0.4}}$ & $\Delta E_{D,K_{0.6}}$ & $\Delta E_{D,G}$ & $\Delta E_{D,GS}$ \\
		\hline
		\multirow{4}{1.5em}{(0,0)} & 18 & 0.00 & 0.00 & 0.00 & 0.00 & 0.00 \\
		& 19 & 0.00 & 0.00 & 0.00 & 0.06 & 0.00 \\
		& 19.9 & 0.22 &	0.86 & 1.86 & 13.24 & 1.40 \\
		& 19.99 & 2.30 & 6.06 & 11.59 & 42.6 & 8.84 \\
		\hline
		\multirow{4}{2em}{(0,0.5)} & 18 & 0.00 & 0.00 & 0.00 & 0.00 & 0.00 \\
		& 19 & 0.00 & 0.00 & 0.00 & 0.03 & 0.00 \\
		& 19.9 & 0.10 & 0.38 & 0.75 & 5.59 & 0.69 \\
		& 19.99 & 1.07 & 2.59 &	4.42 & 15.65 & 4.18 \\
		\hline
		\multirow{4}{2em}{(1.5,0)} & 18 & 0.00 & 0.00 & 0.00 & 0.00 & 0.00 \\
		& 19 & 0.00 & 0.00 & 0.00 & 0.00 & 0.00 \\
		& 19.9 & 0.06 &	0.39 & 1.02 & 10.26 & 0.84 \\
		& 19.99 & 1.80 & 5.14 & 10.11 & 38.14 & 7.84 \\
		\hline
		\multirow{4}{3em}{(1.5,0.5)} & 18 & 0.00 & 0.00 & 0.00 & 0.00 & 0.00 \\
		& 19 & 0.00 & 0.00 & 0.00 & 0.00 & 0.00 \\
		& 19.9 & 0.01 & 0.09 & 0.25 & 3.51 & 0.27 \\
		& 19.99 & 0.69 & 1.93 & 3.49 & 12.87 & 3.38 \\
		\hline
	\end{tabular}
	\label{EfficiencyOptimal2}
\end{table*}

The observations from the new no-show function is also in line with Proposition 3 in Liu (2016). As demand rate increases, restricting scheduling window will bring more efficiency due to delay reduction and in turn controlling no-show rate. In addition, as patient rejection penalty or ancillary task revenue increases, efficiency gain decreases since optimal scheduling window will become larger and more patients will be able to get an appointment within scheduling window. 

\section{Efficiency Gains Resulted from Adopting $K^*$ When Practices May Adjust Panel Size and Overbooking Level}

Liu (2016) proposed the following optimization model to find the optimal value of demand rate (provider panel size), $\lambda$, and daily service rate (level of overbooking), $\mu$:
\begin{equation}
max_{\lambda,\mu} \lambda \sum_{j=0}^{\infty}\Pi_{j}(\lambda,\mu)q{j} + \mu\xi(1-\rho) - \omega(\mu)
\end{equation}
Where $\rho$ is the traffic intensity of the system and, $\Pi_{j}(\lambda,\mu)$ is the steady-state probability of having j patients in the system given demand rate $\lambda$ and service rate $\mu$, and $\omega(\mu)$ is defined as $a\ast[(\mu-M)^+]^2$ in which $M=20$ is the regular daily capacity and $a \in \{0.2,2\}$ is the provider overtime cost parameter.

Once the clinic finds the optimal demand rate and service rate by solving (6), it can further improve efficiency by adopting an optimal scheduling window through solving (3). Table \ref{EfficiencyLevers} shows how much the clinic will gain in terms of efficiency by adopting optimal appointment scheduling window. We also evaluate both queuing systems in terms of service level which is defined as the probability that no more than $K^*$ patients be in the system with infinity scheduling window and optimal panel size and overbooking level.
\begin{table*}[!htbp]
	\caption{Efficiency Gains (\%) Resulted from Adopting $K^*$ with Other Operational Levers}
	\centering
	\begin{tabular}{|c c c c c c c c c c|} 
		\hline
		$(\theta,\xi)$ & No-show & $\Delta E^{(\lambda^*,\mu)}_{M}$ & $\alpha^{(\lambda^*,\mu)}_{M}$ & $\Delta E^{(\lambda^*,\mu)}_{D}$ & $\alpha^{(\lambda^*,\mu)}_{D}$ & $\Delta E^{(\lambda^*,\mu^*)}_{M}$ & $\alpha^{(\lambda^*,\mu^*)}_{M}$ & $\Delta E^{(\lambda^*,\mu^*)}_{D}$ & $\alpha^{(\lambda^*,\mu^*)}_{D}$ \\
		\hline
		\multirow{4}{1.5em}{(0,0)} & $K_{0.2}$ & 0.03 & 0.87 & 0.20 & 0.70 & 0.13 & 0.82 & 0.11 & 0.79 \\
		& $K_{0.4}$ & 0.12 & 0.88 & 0.15 & 0.87 & 0.21 & 0.83 & 0.15 & 0.82 \\
		& $K_{0.6}$ & 0.20 & 0.87 & 0.21 & 0.79 & 0.29 & 0.85 & 0.21 & 0.81 \\
		& $G$ & 0.97 & 0.81 & 0.43 & 0.73 & 0.81 & 0.84 & 0.85 & 0.68 \\
		& $GS$ & 0.29 & 0.80 & 0.06 & 0.92 & 0.15 & 0.87 & 0.20 & 0.86 \\
		\hline
		\multirow{4}{2em}{(0,0.5)} & $K_{0.2}$ & 0.05 & 0.87 & 0.05 & 0.84 & 0.06 & 0.82 & 0.05 & 0.79 \\
		& $K_{0.4}$ & 0.07 & 0.88 & 0.07 & 0.87 & 0.09 & 0.83 & 0.06 & 0.82 \\
		& $K_{0.6}$ & 0.09 & 0.87 & 0.09 & 0.79 & 0.12 & 0.85 & 0.09 & 0.81 \\
		& $G$ & 0.43 & 0.81 & 0.19 & 0.73 & 0.36 & 0.84 & 0.38 & 0.68 \\
		& $GS$ & 0.14 & 0.80 & 0.03 & 0.92 & 0.07 & 0.87 & 0.10 & 0.86 \\
		\hline
		\multirow{4}{2em}{(1.5,0)} & $K_{0.2}$ & 0.01 & 0.99 & 0.01 & 0.97 & 0.02 & 0.98 & 0.02 & 0.97 \\
		& $K_{0.4}$ & 0.01 & 0.99 & 0.02 & 0.98 & 0.02 & 0.98 & 0.02 & 0.98 \\
		& $K_{0.6}$ & 0.01 & 0.99 & 0.02 & 0.98 & 0.02 & 0.99 & 0.02 & 0.99 \\
		& $G$ & 0.13 & 0.96 & 0.06 & 0.98 & 0.09 & 0.97 & 0.16 & 0.99 \\
		& $GS$ & 0.05 & 0.96 & 0.00 & 1.00 & 0.01 & 0.99 & 0.04 & 0.97 \\
		\hline
		\multirow{4}{2em}{(1.5,0.5)} & $K_{0.2}$ & 0.00 & 1.00 & 0.00 & 1.00 & 0.00 & 1.00 & 0.00 & 1.00 \\
		& $K_{0.4}$ & 0.00 & 1.00 & 0.00 & 1.00 & 0.00 & 1.00 & 0.00 & 1.00 \\
		& $K_{0.6}$ & 0.00 & 1.00 & 0.00 & 1.00 & 0.00 & 1.00 & 0.00 & 1.00 \\
		& $G$ & 0.00 & 1.00 & 0.00 & 1.00 & 0.00 & 1.00 & 0.02 & 0.99 \\
		& $GS$ & 0.00 & 0.99 & 0.00 & 1.00 & 0.00 & 1.00 & 0.00 & 0.99 \\
		\hline
	\end{tabular}
	\label{EfficiencyLevers}
\end{table*}

Results in Table \ref{EfficiencyLevers} reveals that as estimated no-show rate increases in function proposed by Kopach et al. (2007), efficiency gain of adopting optimal scheduling window after adjusting panel size and overbooking level will increase too since patients become more sensitive to delay and clinic has to reduce scheduling window.

\section{Efficiency Gains Resulted from All Operational Levers}

The following optimization problem is proposed by Liu (2016) to find the optimal value of scheduling window:
\begin{equation}
max_{\lambda,\mu} \lambda \sum_{j=0}^{\infty}\Pi_{j}(\lambda,\mu,K)q{j} + \mu\xi\Pi_{0}(\lambda,\mu,K) -\lambda\theta\Pi_{K}(\lambda,\mu,K) - \omega(\mu)
\end{equation}
Where $\Pi_{j}(\lambda,\mu)$ is the steady-state probability of having j patients in the system given demand rate $\lambda$ and service rate $\mu$, and $\omega(\mu)$ is defined as $a\ast[(\mu-M)^+]^2$ in which $M=20$ is the regular daily capacity and $a \in \{0.2,2\}$ is the provider overtime cost parameter. The average efficiency gain in M/M/1/K system is almost 0\% and 4\% when $\theta = 0$ and $\theta \neq 0$, respectively. This confirms the conclusion by Liu (2016) that optimizing the appointment scheduling window is a substitute rather than a complement to optimizing operational levers such as the panel size and overbooking level. 

\section{References}
Kopach, R., DeLaurentis, P., Lawley, M., Muthuraman, K., Ozsen, L., Rardin, R., Wan, H., Intrevado, P., Qu, X., Willis, D. 2007. Effects of clinical characteristics on successful open access scheduling. Health care management science, 10(2), 111-124. \\
Liu, N. 2016. Optimal choice for appointment scheduling window under patient no‐show behavior. Production and Operations Management, 25(1), pp.128-142.

\end{document}